\documentclass[11pt]{amsart}

\usepackage{amssymb,amsmath,amsxtra,amscd}
\usepackage[mathscr]{eucal} 

\theoremstyle{plain}
\numberwithin{equation}{section}
\newtheorem{theorem}{Theorem}
\newtheorem{proposition}{Proposition}

\numberwithin{lemma}{section}
\numberwithin{theorem}{section}
\numberwithin{corollary}{section}
\numberwithin{proposition}{section}

\theoremstyle{definition}
\newtheorem{definition}{Definition}
\numberwithin{definition}{section}

\theoremstyle{remark}

\numberwithin{example}{section}

\newcommand{\nc}{\newcommand}
\nc{\mc}{\mathcal}
\nc{\MSV}[1]{\ensuremath{ \Omega^{\on{ch}}_{#1 } }}
\nc{\MSVg}[1]{\ensuremath{ \Omega^{\on{ch},g}_{#1 } }}
\nc{\on}{\operatorname}
\nc{\BunG}{\on{BunG}}
\nc{\BunGH}{\on{\BunG^{H}}}
\nc{\Mg}{\on{M}^{g}}
\nc{\M}{\on{M}}
\renewcommand{\H}{\mc{H}}
\nc{\Cl}{Cl}
\nc{\drva}{\widehat{\Omega}}
\nc{\spec}{\on{Spec}}
\nc{\AutO}{\on{Aut} \mc{O}}
\nc{\vac}{|0\rangle}
\nc{\Z}{\mathbb{Z}}
\nc{\zf}[1]{z^{\frac{1}{#1}}}
\nc{\wf}[1]{w^{\frac{1}{#1}}}
\nc{\Mt}{M^{\sigma}}
\renewcommand{\exp}{\on{exp}}
\nc{\gr}{\on{gr}}
\renewcommand{\P}{\mathcal{P}}

\nc{\al}{\alpha}
\nc{\ol}{\overline}

\begin{document}

\title{Chiral De Rham Complex and Orbifolds}

\author{Edward Frenkel}\thanks{Partially supported by grants from the
Packard Foundation and NSF}
\address{Department of Mathematics  
         University of California, Berkeley, CA 94720}
\email{frenkel@math.berkeley.edu}
\author{Matthew Szczesny}
\address{Department of Mathematics  
         University of Pennsylvania, Philadelphia, PA 19104 }
\email{szczesny@math.upenn.edu}

\date{July 2003}

\begin{abstract}
Suppose that a finite group $G$ acts on a smooth complex variety
$X$. Then this action lifts to the Chiral de Rham Complex $\MSV{X}$ of
$X$ and to its cohomology by automorphisms of the vertex algebra
structure. We define twisted sectors for $\MSV{X}$ (and their
cohomologies) as sheaves of twisted vertex algebra modules supported
on the components of the fixed-point sets $X^{g}, g \in G$. Each
twisted sector sheaf carries a BRST differential and is
quasi-isomorphic to the de Rham complex of $X^{g}$. Putting the
twisted sectors together with the vacuum sector and taking
$G$--invariants, we recover the additive and graded structures of
Chen-Ruan orbifold cohomology. Finally, we show that the orbifold
elliptic genus is the partition function of the direct sum of the
cohomologies of the twisted sectors.
\end{abstract}

\maketitle 

\section{Introduction}

In \cite{MSV} the authors constructed a sheaf of vertex superalgebras
on any smooth complex variety $X$, called the chiral de Rham complex,
and denoted by $\Omega^{\on{ch}}_{X}$ (a chiral algebra version of
this construction has been given in \cite{BD}). The cohomology of
$\Omega^{\on{ch}}_{X}$, $H^{*}(X,\Omega^{\on{ch}}_{X})$, is also a
vertex superalgebra. Thus, one attaches a vertex superalgebra to an
arbitrary smooth algebraic manifold, and when $X$ is Calabi-Yau, it is
believed to be related to the chiral algebra of the topological twist
of the type $A$ sigma-model with target $X$ (at least in the large
radius limit). This explains the recent interest in the chiral de Rham
complex: using it one may hopefully obtain a rigorous mathematical
construction of a conformal field theory associated to $X$ (such as
the sigma model of $X$) without resorting to non-rigorous methods like
the path integral. The chiral de Rham complex also has a number of
interesting applications in various fields of geometry and
representation theory, such as mirror symmetry \cite{Bo}, elliptic
genera \cite{BL1}, \cite{BL2}, as well as others \cite{GMS},
\cite{GMS2}, \cite{GMS3}, \cite{S}. We hope that in the future it may
be used to access other aspects of ``stringy'' geometry.

The aim of this paper is to extend the construction of
$\Omega^{\on{ch}}_{X}$ to orbifolds. Suppose that $X$ carries the
action of a finite group $G$. In physics terms, we would then like to
construct an orbifold conformal field theory (such as the sigma model)
corresponding to the stack $[X/G]$. In the language of vertex
algebras, we must first take the $G$--invariant part of
$H^{*}(X,\Omega^{\on{ch}}_{X})$ and then include the so-called twisted
sectors, one for each conjugacy class in $G$. Let $T$ denote the set
of conjugacy classes of $G$, and $[g]$ the conjugacy class of an
element $g \in G$. In this paper we give a geometric construction of
these twisted sectors as follows. For each $g \in G$, let $X^{g}$
denote the fixed-point set of $g$ on $X$. We construct a twisted
version of $\Omega^{\on{ch}}_{X}$, denoted $\Omega^{\on{ch},g}_{X}$,
as a sheaf supported on $X^{g}$. We show that $\Omega^{\on{ch},g}_{X}$
is a sheaf of $g$--twisted $\Omega^{\on{ch}}_{X}$--modules, and that
furthermore, $H^{*}(X,\Omega^{\on{ch},g}_{X})$ is a $g$--twisted
module for $H^{*}(X,\Omega^{\on{ch}}_{X})$. Taking $G$--invariants, we
show that $H^{*}(X,\Omega^{\on{ch},g}_{X})^{C(g)}$ is an ordinary
vertex algebra module for $H^{*}(X,\Omega^{\on{ch}}_{X})^{G}$, where
$C(g)$ denotes the centralizer of $g$. Furthermore, if $[g] = [g']$,
then we have a canonical isomorphism
\[
H^{*}(X,\Omega^{\on{ch},g'}_{X})^{C(g')} \simeq
H^{*}(X,\Omega^{\on{ch},g}_{X})^{C(g)}
\]
of $H^{*}(X,\Omega^{\on{ch}}_{X})^{G}$--modules. In particular,
$\Omega^{\on{ch},1}_{X} = \MSV{X}$. The space
\[
 \mathscr{H} = \bigoplus_{[g] \in T} H^{*}(X,\Omega^{\on{ch},g}_{X})^{C(g)}
\]
is then the conjectural space of states of the orbifold model. If $X$
is Calabi-Yau and $G$ acts preserving the volume form, then the
untwisted sector of $\mathcal{H}$ has the structure of an $N=2$
topological vertex algebra.

As shown in \cite{MSV}, $\Omega^{\on{ch}}_{X}$ carries a BRST differential,
denoted $d^{\on{ch}}$, and there is an inclusion of the usual de Rham
complex
\[
	(\Omega_{X}, d_{DR}) \hookrightarrow (\Omega^{\on{ch}}_{X},
	d^{\on{ch}})
\]
which is a quasi-isomorphism. In the orbifold setting, each sheaf
$\Omega^{\on{ch},g}_{X}$ carries a BRST differential $d^{\on{ch},g}$, and
there is a quasi-isomorphism
\[
	(\Omega_{X^{g}}, d_{DR}) \hookrightarrow (\Omega^{\on{ch},g}_{X},
	d^{\on{ch},g}).
\]
It follows that the hypercohomology of $\Omega^{\on{ch},g}$ is
isomorphic to
\[
	{\mathbb H}^{*}(X,\Omega^{\on{ch},g}_{X})^{C(g)} \cong
	H^{*}(X^{g}, \mathbb{C})^{C(g)} \cong
	H^{*}(X^{g}/C(g), \mathbb{C})
\]
Combining all of the twisted sectors together, we see that
\begin{equation}    \label{chen-ruan}
	\bigoplus_{[g] \in T} {\mathbb
	H}^{*}(X,\Omega^{\on{ch},g}_{X})^{C(g)} \cong \bigoplus_{[g] \in
	T} H^{*}(X^{g}/C(g), \mathbb{C})
\end{equation}
is thus isomorphic to the Chen-Ruan orbifold cohomology of $[X/G]$
(see \cite{CR}). The important point of our construction is that ${\mc
H}$, and hence the right hand side of \eqref{chen-ruan} acquires a
natural gradation from the vertex algebra point of view. We show that
this gradation coincides with the somewhat mysterious gradation
(shifted by ``age'' or ``fermionic shift'') on the orbifold cohomology
introduced by Chen and Ruan. Thus, we obtain an isomorphism between
the left hand side of \eqref{chen-ruan} and the orbifold cohomology at
the level of graded vector spaces. We expect that the product
structure on the orbifold cohomology may also be constructed by the
vertex algebra methods, along the lines of \cite{MS} (where the
question of constructing of the product structure on the ordinary de
Rham cohomology in terms of the chiral de Rham complex was
considered). We hope to return to this question in a future
publication.

In \cite{GMS},\cite{GMS2}, a purely even version of the chiral de Rham
complex, the sheaf of chiral differential operators, was defined for
those varieties $X$ on which the second Chern character vanishes (see
\cite{BD} for a chiral algebra version). If a finite group $G$ acts on
such a variety, then we can define in a similar way sheaves of twisted
modules over the sheaf of chiral differential operators.

It has been shown in \cite{BL} that the elliptic genus of $X$ may be
expressed as the trace of the operator $y^{J_{0}} q^{L_{0}}$ on
$H^{*}(X,\Omega^{\on{ch}}_{X})$. Here we prove an analogue of this
statement in the orbifold setting. Namely, we show that the orbifold
elliptic genus of $[X/G]$ may be expressed as the trace of $y^{J_{0}}
q^{L_{0}} $ on $\mathcal{H}$. The latter may be viewed as an invariant
of the corresponding orbifold conformal field theory.

Our results can be generalized to more general orbifolds, which are not
global quotients.

\medskip

\noindent{\bf Acknowledgments} We thank D. Kaledin for suggesting to
us the possibility that orbifold cohomologies could be obtained using
the chiral de Rham complex. We are also grateful to F. Malikov for
useful comments.

The results of this paper were reported at the Workshop ``Integration
on Arc Spaces, Elliptic Genus and Chiral de Rham Complex'' in Banff in
June of 2003. We thank the organizers for the invitation. At this
conference, A. Vaintrob informed us that he has also constructed an
orbifold version of the chiral de Rham complex, and obtained results
similar to ours.

\section{Vertex Algebras and Twisted Modules}

In this paper we will use the language of vertex superalgebras, their
modules, and twisted modules. For an introduction to vertex algebras
and their modules \cite{FLM,KAC,FBZ}, and for background on twisted
modules, see \cite{FFR,Dong,DLM,FS}.

We recall that a conformal vertex superalgebra is a
$\mathbb{Z}_+$--graded super vector space $$V = \bigoplus_{n=0}^\infty
V_n,$$ $$V_{n}=V_{n}^{\ol{0}} \oplus V_{n}^{\ol{1}}$$ together with a
vacuum vector $\vac \in V_0^{\ol{0}}$, an even translation operator $T$
of degree $1$, a conformal vector $\omega \in V_2^{\ol{0}}$ and an even
linear map

\begin{align*}
Y: V &\to \on{End} V[[z^{\pm 1}]], \\ A &\mapsto Y(A,z) = \sum_{n \in
\Z} A_{(n)} z^{-n-1}.
\end{align*}
These data must satisfy certain axioms
(see \cite{FLM,KAC,FBZ}). In what follows we will denote the
collection of such data simply by $V$, and the parity of an element $A
\in V$ homogeneous with respect to the $\mathbb{Z}/2 \mathbb{Z}$ grading by
$p(A)$.

A vector superspace $M$ is called a $V$--module if it is equipped
with an even linear map
\begin{align*}
Y^M: V &\to \on{End} M[[z^{\pm 1}]], \\
A &\mapsto Y^M(A,z) = \sum_{n \in \Z} A^M_{(n)} z^{-n-1}
\end{align*}
such that for any $v \in M$ we have $A^M_{(n)} v = 0$ for large enough
$n$. This operation must satisfy the following axioms:

\begin{itemize}

\item $Y^M(\vac,z) = \on{Id}_M$;

\item For any $v \in M$ and homogeneous $A,B \in V$ there exists an
element $$f_v \in M[[z,w]][z^{-1},w^{-1},(z-w)^{-1}]$$ such that the
formal power series
$$Y^M(A,z) Y^M(B,w) v, \quad (-1)^{p(A)p(B)}Y^{M}(B,w) Y^{M}(A,z)v ,
\qquad \on{and} $$ $$ \qquad Y_M(Y(A,z-w) B,w) v$$ are expansions of $f_v$
in $M((z))((w))$, $M((w))((z))$ and $M((w))((z-w))$, respectively.

\end{itemize}

The power series $Y^M(A,z)$ are called vertex operators. We write the
vertex operator corresponding to $\omega$ as
\[
	Y^M(\omega,z) = \sum_{n \in \mathbb{Z}} L^M_{n} z^{-n-2},
\]
where $L^M_n$ are linear operators on $V$ generating the Virasoro
algebra. Following \cite{Dong}, we call $M$ \emph{admissible} if
$L^{M}_{0}$ acts semi-simply with integral eigenvalues.

Now let $\sigma_{V}$ be a conformal automorphism of $V$, i.e., an even
automorphism of the underlying vector superspace preserving all of the
above structures (so in particular $\sigma_{V}(\omega) = \omega$). We
will assume that $\sigma_{V}$ has finite order $m>1$. A vector space
$M^\sigma$ is called a $\sigma_V$--{\em twisted} $V$--module (or
simply twisted module) if it is equipped with an even linear map
\begin{align*}
Y^{M^\sigma}: V &\to \on{End} M^\sigma[[z^{\pm \frac{1}{m}}]], \\ A
&\mapsto Y^{M^{\sigma}}(A,\zf{m}) = \sum_{n \in \frac{1}{m}\Z}
A^{M^\sigma}_{(n)} z^{-n-1}
\end{align*}
such that for any $v \in M^{\sigma}$ we have $A^{M^\sigma}_{(n)} v =
0$ for large enough $n$. Please note that we use the notation
$Y^{M^{\sigma}}(A,\zf{m})$ rather than $Y^{M^{\sigma}}(A,z)$ in the
twisted setting. This operation must satisfy the following axioms (see
\cite{FFR,Dong,DLM,LI,FS}):

\begin{itemize}

\item $Y^{M^{\sigma}}(\vac,\zf{m}) = \on{Id}_{\Mt}$;

\item For any $v \in \Mt$ and homogeneous $A,B \in V$, there exists an
element
\[
f_{v} \in \Mt [[\zf{m}, \wf{m} ]][z^{- \frac{1}{m}}, w^{-
\frac{1}{m}},(z-w)^{-1}]
\]
such that the formal power series $$Y^{\Mt}(A,\zf{m}) Y^{\Mt}(B,
\wf{m})v, \quad (-1)^{p(A)p(B)} Y^{\Mt}(B,\wf{m}) Y^{\Mt}(A, \zf{m})v,
\qquad \on{and} $$ $$ \qquad Y^{\Mt}(Y(A,z-w) B, \wf{m})v$$ are expansions
of $f_{v}$ in $\Mt((\zf{m}))((\wf{m}))$, $\Mt((\wf{m}))((\zf{m})) $
and $\Mt((\wf{m}))((z-w))$, respectively.

\item If $A \in V$ is such that $\sigma_V(A) = e^{\frac{2\pi i k}{m}}
A$, then $A^{M^\sigma}_{(n)} = 0$ unless $n \in \frac{k}{m} + \Z$.

\end{itemize}

The series $Y^{M^\sigma}(A,z)$ are called twisted vertex operators.
In particular, the Fourier coefficients of the twisted vertex operator
\[
	Y^{M^\sigma}(\omega,\zf{m}) = \sum_{n \in \mathbb{Z}}
	L^{M^\sigma}_{n} z^{-n-2},
\]
generate an action of the Virasoro algebra on $M^\sigma$. 

\section{Chiral de Rham Complex}

In this section we review the construction of the chiral deRham
complex due to \cite{MSV}. We follow the presentation of \cite{FBZ},
Sect.~17.4.

\subsection{The vertex superalgebra $\drva_{N}$}

Let $\H_{N}$ denote the infinite-dimensional Heisenberg algebra with
generators $a^{i}_{n}, \; b^{j}_{m}, \quad i, j = 1, \cdots, N$ $n, m
\in \mathbb{Z}$, and relations
\[
	[a^{i}_{n}, b^{j}_{m}] = \delta_{i, j} \delta_{n, -m}.
\]
Let $\H^{+}_{N}$ denote the abelian Lie subalgebra generated by
$a^{i}_{n}, \; n \geq 0$,  $ b^{j}_{m}, \; m > 0$. Let
\[
	\on{V}_{N} = \on{Ind}^{\H_{N}}_{\H^{+}_{N}} \mathbb{C},
\]
where $\mathbb{C}$ denotes the trivial representation of $\H^{+}_{N}$.
Then $\on{V}_{N}$ has the structure of a vertex algebra
(see \cite{FBZ}). It is shown in \cite{MSV} that the vertex algebra
structure on $\on{V}_{N}$ can be naturally extended to its completion
\[
	\widehat{\on{V}}_{N} = \mathbb{C}[[b^{i}_{0}]]
	\underset{\mathbb{C}[b^{i}_{0}]}\otimes \on{V}_{N}.
\]

Let $\Cl_{N}$ be the infinite-dimensional Clifford algebra with
generators $\psi^{i}_{n}, \; \phi^{j}_{m}, \quad i, j = 1,\cdots,N$,
$n, m \in \mathbf{Z}$, and relations
\[
	[\psi^{i}_{n}, \phi^{j}_{m}]_{+} = \delta_{i,j} \delta_{n, -m}
\]
Let $\Cl^{+}_{N}$ denote the abelian Lie sub-superalgebra spanned by
$\psi^{i}_{n}, \; n \geq 0$, $ \phi^{j}_{m}, \; m > 0$. Let
\[
	\Lambda_{N} = \on{Ind}^{\Cl_{N}}_{\Cl^{+}_{N}} \mathbb{C},
\]
where $\mathbb{C}$ denotes the trivial representation of
$\Cl^{+}_{N}$.  Then $\Lambda_{N}$ has the structure of a vertex
superalgebra (see, e.g., \cite{FBZ}).

The tensor product of vertex (super) algebras has the structure of
a vertex (super) algebra in a natural way. Let
\[
	\drva_{N} = \widehat{\on{V}}_{N} \otimes \Lambda_{N}
\]
Let us denote by $Y: \drva_{N} \mapsto \on{End}(\drva_{N})[[z,z^{-1}]]$
the corresponding vertex operator.

\subsection{Superconformal structure} \label{SCF}

Define the vectors $L,\; J, \; Q, \; G \; \in \drva_{N} $ as follows:
\[
	L = \sum^{N}_{i=1} (b^{i}_{-1} a^{i}_{-1} + \phi^{i}_{-1}
	\psi^{i}_{-1}),
\]
\[
	J = \sum^{N}_{i=1} \phi^{i}_{0} \psi^{i}_{-1}, \qquad Q=
	\sum^{N}_{i=1} a^{i}_{-1} \phi^{i}_{0}, \qquad G =
	\sum^{N}_{i=1} \psi^{i}_{-1} b^{i}_{-1}
\]
The corresponding fields, $L(z) = Y(L,z), J(z) = Y(J,z), Q(z)= Y(Q,z),
G(z)= Y(G,z)$ are:
\[
	L(z) = \sum^{N}_{i=1} (:\partial_{z} b^{i}(z) a^{i}(z) : +
	:\partial_{z} \phi^{i}(z) \psi^{i}(z):),
\]
\[
	J(z) = \sum^{N}_{i=1} :\phi^{i}(z) \psi^{i}(z) : , \; Q(z) =
	\sum^{N}_{i=1} :a^{i}(z) \phi^{i}(z) : , \; G(z) =
	\sum^{N}_{i=1} :\psi^{i}(z) \partial_{z} b^{i}(z) :
\]
It is known from \cite{FMS} (see also \cite{MSV}) that the fields
$L(z), J(z), Q(z), G(z)$ have the OPE's of an $N=2$ superconformal
algebra of central charge $0$. These OPE's are:

$$
L(z)L(w)=\frac{2L(w)}{(z-w)^2}+\frac{\partial_{w} L(w)}{z-w}
$$
$$  
J(z)J(w)=\frac{N}{(z-w)^2};\ 
L(z)J(w)=-\frac{N}{(z-w)^3}+\frac{J(w)}{(z-w)^2}+\frac{\partial_{w}J(w)}{z-w}
$$
$$
G(z)G(w)=0;\ L(z)G(w)=\frac{2G(w)}{(z-w)^2}+\frac{\partial_{w}G(w)}{z-w};\ 
J(z)G(w)=-\frac{G(w)}{z-w}
$$
$$
Q(z)Q(w)=0;\ L(z)Q(w)=\frac{Q(w)}{(z-w)^2}+\frac{\partial_{w}Q(w)}{z-w};\ 
J(z)Q(w)=\frac{Q(w)}{z-w}
$$
$$
Q(z)G(w)=\frac{N}{(z-w)^3}+\frac{J(w)}{(z-w)^2}+\frac{L(w)}{z-w}
$$

\subsection{Charges and gradings}

The operator 
\[
	L_{0} = \sum^{N}_{i=1} \sum_{n \in \mathbb{Z}} (:a^{i}_{n}
	b^{i}_{-n}: + :\phi^{i}_{n} \psi^{i}_{-n}:)
\]
defines a $\mathbb{Z}_{\geq 0}$ grading on $\drva_{N}$ by conformal
weights. Let
\[
\drva_{N,k} = \{v \in \drva_{N} \vert L_{0} v = kv \}
\]
The operator
\[
	F = J_{0} = \sum^{N}_{i=1} \sum_{n \in \mathbb{Z}}
	:\phi^{i}_{n} \psi^{i}_{-n}:
\]
defines a $\mathbb{Z}$ grading on $\drva_{N}$ by fermionic
charge. Let 
\[
	\drva^{p}_{N} = \{ v \in \drva_{N} \vert F v = p v \}
\]
$L_{0}$ and $J_{0}$ commute, and thus define a $\mathbb{Z}_{\geq 0}
\times \mathbb{Z}$ bigrading on $\drva_{N}$. We use the notation
$\drva^{p}_{N,k}$ for the bigraded components, and call it the
component of bidegree $(k,p)$. Let 
\[
	d^{\on{ch}} = - Q_{0} = - \sum^{N}_{i=1} \sum_{n \in \mathbb{Z}}
	a^{i}_{n} \phi^{i}_{-n}
\]
be the BRST or deRham differential of $\drva_{N}$. It has bidegree
$(0,1)$ and $d^{\on{ch}} \circ d^{\on{ch}} = 0$.

\subsection{Construction of $\MSV{X}$} \label{constructMSV}

Let $D_{N}$ be the formal $N$--dimensional disk $$D_N = \spec
\mathbb{C}[[t_{1}, \cdots, t_{N}]].$$ Let $W_{N}$ denote the
topological Lie algebra of vector fields on $D_{N}$, $W^{+}_{N}$ the
Lie subalgebra of vector fields vanishing at the origin, and
$\AutO_{N}$ the unipotent proalgebraic group of continuous algebra
automorphisms of $\mathbb{C}[[t_{1}, \cdots, t_{N}]]$. We have that
$\on{Lie}(\AutO_{N}) = W^{+}_{N}$, and $(W_{N}, \AutO_{N})$ forms a
Harish-Chandra pair.

Let $U(\drva_{N})$ denote the Lie algebra of Fourier coefficients of
vertex operators on $\drva_{N}$. It follows from the results in
\cite{MSV} that there exists a Lie algebra homomorphism
\[
	\nu: W_{N} \mapsto U(\drva_{N})
\]
defined by the formula
\begin{equation} \label{Wrep}
	f(t_{i}) \partial_{t_{j}} \mapsto \on{Res}_{z} \left( :f(b^{i}(z))
	a^{j}(z): + \sum^{N}_{k=1} : (\partial_{t_{k}} f)(b^{i}(z))
	\phi^{k}(z) \psi^{j}(z) : \right)
\end{equation}
Since $U(\drva_{N})$ acts on $\drva_{N}$, we obtain a representation
of $W_{N}$ on $\drva_{N}$ by residues of vertex operators. Over
$W^{+}_{N}$ this action can be exponentiated to an action of
$\AutO_{N}$ by vertex algebra automorphisms. 

The infinitesimal action of $\AutO_{N}$ on the generating fields of
$\drva_{N}$ is given by
\begin{equation} \label{action1}
	[\nu(f(t_{i}) \partial_{t_{j}} ), a^{k}(w)] =
	:(\partial_{t_{k}} f)(b^{i}(w))a^{j}(w): + \sum^{N}_{l=1}
	:(\partial_{t_{l}} \partial_{t_{k}} f)(b^{i}(w)) \phi^{l}(w)
	\psi^{j}(w):
\end{equation}
\begin{equation} \label{action2}
	[\nu(f(t_{i}) \partial_{t_{j}} ), b^{k}(w)] = \delta_{j,k}
	f(b^{i}(w))
\end{equation}
\begin{equation} \label{action3}
	[\nu(f(t_{i}) \partial_{t_{j}} ), \phi^{k}(w)] = - \delta_{j,k}
	\sum^{N}_{l=1} (\partial_{t_{l}} f)(b^{i}(w)) \phi^{l}(w)
\end{equation}
\begin{equation} \label{action4}
	[\nu(f(t_{i}) \partial_{t_{j}} ), \psi^{k}(w)] =
	(\partial_{t_{k}} f)(b^{i}(w)) \psi^{j}(w)
\end{equation}

Let $X$ be a smooth complex variety of dimension $N$. Then $X$ carries
a canonical Harish-Chandra structure $\widehat{X}$
\[
	\widehat{X} = \{ (x, (t_{i,x})_{i=1 \cdots N}) \}
\]
where $x \in X$ and $(t_{i,x})_{i=1,\cdots,N}$ is a set of topological
generators of the maximal ideal of $\widehat{\mc{O}}_{x}$. It is clear
that $\widehat{X}$ is an $\AutO_{N}$--principal bundle over $X$ via
the map $(x, (t_{i,x})) \mapsto x$ (see \cite{FBZ}, where it is
denoted by ${\mathcal A}ut_X$). Moreover, the Lie algebra $W_N$ acts
on $\widehat{X}$ simply transitively, extending the fiberwise action of
$W^+_N$ induced by the action of $\AutO_N$. Let
\[
	\widetilde{\Omega}^{\on{ch}}_{X} = \widehat{X}
	\underset{\AutO_{N}} \times \drva_{N}
\]
Thus, the fiber of $\widetilde{\Omega}^{\on{ch}}_{X}$ at $x$ is the
quotient of the space of pairs
$(t_{i,x}, v)$, where $t_{i,x}$ are formal coordinates on $\spec
\widehat{\mc{O}}_{x}$, $v \in \drva_{N}$, by the equivalence
relation
\[
	(t_{i,x}, v) \sim (\rho(t_{i,x}), \rho(v)), \quad \rho \in
	\AutO_{N}.
\]
As explained in \cite{FBZ}, the sheaf of sections of
$\widetilde{\Omega}^{\on{ch}}_{X}$ is an ${\mathcal O}_X$--module
equipped with a flat connection, i.e., it is a
$\mc{D}_{X}$--module. The connection is constructed from the action of
the translation vector fields $\partial_{t_{i}}$ on $\drva_{N}$.  The
chiral de Rham complex $\MSV{X}$ is by definition the sheaf of
horizontal sections of $\widetilde{\Omega}^{\on{ch}}_{X}$. It is a
sheaf of vertex superalgebras whose cohomology also possesses a vertex
superalgebra structure, as explained in \cite{MSV}, \cite{Bo}.

We proceed to give an explicit description of the $\mc{D}_{X}$--module
structure of $\widetilde{\Omega}^{\on{ch}}_{X}$ in local
coordinates. Let $x \in X$, and choose formal coordinates $t_{i,x}$ at
$x$. Denote $\on{Spec} \widehat{\mc{O}}_{x}$ by $D_{x}$. Then we
obtain a trivialization $\tau$ of the restriction of the bundle
$\widetilde{\Omega}^{\on{ch}}_{X}$ to $D_{x}$. It follows from the
formula \eqref{Wrep} that the connection operator $\nabla_{t_{i,x}}$
is given by
$$
\partial_{t_{i,x}} - a^{i}_{0}
$$
Using the trivialization $\tau$ we can therefore identify the space of
horizontal sections of $\widetilde{\Omega}^{\on{ch}}_{X}$ on $D_{x}$ with
those elements $f(t_{i,x}, a^{i}_{n},
b^{i}_{n},\psi^{i}_{n},\phi^{i}_{n})$ of
$$
\mathbb{C}[[t_{1,x}, \cdots, t_{N,x}]] \widehat{\otimes} \drva_{N}
$$
which are ``functions of $t_{i,x} + b^{i}_{0}$'', i.e. those $f$
having the property that they are equal to $0$ after the formal
substitution $t_{i,x} = - b^{i}_{0}, \quad i= 1,\cdots,N$. Evaluation
at $x$, i.e., the map
\[
	\on{ev}_{x}:
	\Gamma_\nabla(D_{x},\widetilde{\Omega}^{\on{ch}}_{X}) \mapsto
	(\widetilde{\Omega}^{\on{ch}}_{X})_{x},
\]
\[
	f \mapsto f \vert_{t_{i,x}=0}
\]
induces an isomorphism between the space of horizontal sections of
$\widetilde{\Omega}^{\on{ch}}_{X}$ on $D_{x}$
and the fiber of $\widetilde{\Omega}^{\on{ch}}_{X}$ at $x$.

The sheaf $\MSV{X}$ is not an
  $\mathcal{O}_X$--module, but it carries a non-associative action of
  $\mc{O}_{X}$
\begin{equation} \label{Oxaction}
	\mu: \mc{O}_{X} \times \MSV{X} \rightarrow
	\MSV{X}
\end{equation}
defined as follows. For $f(t_{i,x}) \in \mathbb{C}[[t_{1,x}, \cdots, t_{N,x}
]]$, denote by $f_{n}$ the operator
\[
\on{Res}_{z} z^{-n-1}
Y(f(b^{i}_{0}),z) \in \on{End} \drva_{N}.
\]
Let $\widetilde{f}_{n} \in \on{End} \drva_{N}[[t_{i,x}]]$ be the
operator obtained from $f_{n}$ by formally replacing $b^{i}_{0}$ by
$b^{i}_{0} + t_{i,x}$.  The action of $f$ on $\MSV{X}(D_{x}) =
\Gamma_\nabla(D_{x},\widetilde{\Omega}^{\on{ch}}_{X}) $ is given by
multiplication by $\widetilde{f}_{0}$. Alternatively, using the
isomorphism $\on{ev}_{x}$, it is given by acting on
$(\widetilde{\Omega}^{\on{ch}}_{X})_{x}$ by $f_{0}$. This action
(\ref{Oxaction}) is non-associative since $f_{0} g_{0} \ne (fg)_{0}$.

\subsection{Filtration on $\MSV{X}$} \label{filtrationMSV}

In \cite{MSV}, the authors introduce a filtration on $\MSV{X}$, which
we now recall. We begin by defining at filtration on
$\drva_{N}$. Note that $\drva_{N}$ is spanned by monomials in the letters
$a^{i}_{n}, \psi^{i}_{n} \,$ $n < 0$ and $b^{i}_{m}, \phi^{i}_{m}, m
\leq 0$. We will call these letters our \emph{alphabet}. We define a
partial order on monomials as follows.
\begin{itemize}
	\item Define a partial order on the alphabet by the rules: $a
	> b, a > \psi, a > b, \psi > \phi , \psi > b, \phi > b$ and
	$x^{i}_{n} > x^{j}_{m}$, x being $a, b, \phi$ or $\psi$
	
	\item Extend the partial order $>$ to all monomials
	lexicographically such that $b^{i}_{0}$ is equivalent to the
	``empty'' letter, i.e. all monomials that differ by products
	of $b^{i}_{0}$ are considered equivalent.
\end{itemize}

This partial order uniquely determines an exhaustive increasing filtration
\[
	F^{0} \drva_{N} \subset F^{1} \drva_{N} \subset F^{2}
	\drva_{N} \cdots
\]
such that $F^i \neq F^{i+1}$ by using the following two rules (here,
denote by $\gr_{F} \MSV{X}$ the associated graded object and $\pi:
\drva_{N} \mapsto \gr_{F} \MSV{X}$ the projection map):

\begin{itemize}
	\item If $M_{1}$ and $M_{2}$ are two monomials, and $M_{1} <
	M_{2}$, then $\pi(M_{1})$ lies in a strictly lower component
	in $\gr_{F} \MSV{X}$ than $\pi(M_{2})$.  \item If neither $M_{1} <
	M_{2}$ nor $M_{1} > M_{2}$ then $\pi(M_{1})$ and $\pi(M_{2})$
	lie in the same graded component.
\end{itemize}

It follows from the formulas \eqref{action1}--\eqref{action4} that this
filtration is $\AutO_{N}$--stable, and therefore induces a filtration
on $\MSV{X}$.

The action $\mu$ (see formula \eqref{Oxaction}) induces an honest
$\mc{O}_{X}$--module structure on $\gr_{F} \MSV{X}$. Furthermore, if
we fix a conformal weight $k \geq 0$, then the sheaves
\begin{equation} \label{assgraded}
	F^{l} \MSV{X} \cap \Omega^{\on{ch}}_{X,k} / F^{l-1} \MSV{X} \cap
	\Omega^{\on{ch}}_{X,k}
\end{equation}
are coherent $\mc{O}_X$--modules isomorphic to tensor products of
exterior and symmetric powers of the tangent and cotangent bundles.

We can give a precise description of $\gr_{F} \MSV{X}$ by using a
generating function, following the notation in \cite{BL}. Note that
$L_{0}, J_{0}$ are $\AutO_{N}$--invariant operators on $\drva_{N}$ and
therefore give rise to global sections of $ \mc{E}nd \MSV{X}$ which
act semi-simply, and which we also denote $L_{0}$ and
$J_{0}$. Furthermore, $L_{0}$ and $J_{0}$ preserve the filtration
$F^{l} \MSV{X}$, and therefore act on $\gr_{F} \MSV{X}$, which splits
into eigenbundles of these operators. We will write
\[
F^{l} \MSV{X} \cap \Omega^{\on{ch}}_{X,k} / F^{l-1} \MSV{X} \cap
	\Omega^{\on{ch}}_{X,k} = \sum_{i} q^{r_{i}} y^{s_{i}}
	W_{r_{i}, s_{i}}
\]
if \eqref{assgraded} is isomorphic to $\oplus_{i} W_{r_{i},s_{i}}$,
and if $L_{0}$ (resp. $J_{0}$) acts on $W_{r_{i},s_{i}}$ by $r_{i}$
(resp. $s_{i}$). If $t$ is a formal variable and $W$ is a vector
bundle, introduce the notation
\[
	\Lambda^{\bullet}_{t} W = \sum_{k \geq 0} t^{k} \Lambda^{k} W 	
\]
and similarly
\[
	\on{Sym}^{\bullet}_{t} W = \sum_{k \geq 0} t^{k} \on{Sym}^{k}
	W.
\]
We then have the following equality (see \cite{BL2})

\begin{proposition}
\begin{align*}
	\gr_{F} \MSV{X} = & \bigoplus_{k,l} F^{l} \MSV{X} \cap
	\Omega^{\on{ch}}_{X,k} / F^{l-1} \MSV{X} \cap
	\Omega^{\on{ch}}_{X,k} \\
			& = \bigotimes_{k \geq 1} \left( 
	\Lambda^{\bullet}_{yq^{k-1}} \Omega_{X} \otimes
	\Lambda^{\bullet}_{y^{-1}q^{k}} T_{X} \otimes
	\on{Sym}^{\bullet}_{q^{k}} \Omega_{X} \otimes
	\on{Sym}^{\bullet}_{q^{k}} T_{X} \right),
\end{align*}
where $T_{X}$ and $\Omega_{X}$ are the holomorphic tangent and
cotangent bundles of $X$ respectively.
\end{proposition}

\section{Twisted Sectors for $\MSV{X}$}

Suppose that $X$ carries an action of a finite group $G$. Then this
action extends naturally to an action on $\widehat{X}$ commuting with
the actions of $\AutO_{N}$ and $W_N$, and therefore gives rise to a
$G$--equivariant structure on $\widetilde{\Omega}^{\on{ch}}_{X}$. For
every $g \in G$ there exists a canonical isomorphism $S_{g} :
\widetilde{\Omega}^{\on{ch}}_{X} \mapsto g^{*}
\widetilde{\Omega}^{\on{ch}}_{X}$. Fibrewise, it is given by the
formula
\begin{align*}
	S_{g}: (\widetilde{\Omega}^{\on{ch}}_{X})_{p} & \mapsto
		  (\widetilde{\Omega}^{\on{ch}}_{X})_{g(p)} \\ (x_{i}
		 , v) & \mapsto (x_i \circ g^{-1}, v).
\end{align*}
These isomorphisms satisfy $S_{h} \circ S_{g} = S_{hg}$

This action of $G$ on $\widetilde{\Omega}^{\on{ch}}_{X}$ takes
horizontal sections to horizontal sections, inducing an action of $G$
on the sheaf $\MSV{X}$. This action respects the vertex algebra
structure, and thus we obtain that $G$ acts on the sheaf $\MSV{X}$ by
vertex algebra automorphisms. This action induces an action of $G$ on
$H^{*}(X, \MSV{X})$ by vertex algebra automorphisms. We also note that
this action is compatible with all gradings as well as the filtration
$F^{l} \MSV{X}$.

In this section we construct on each component $X^g_\al$ of the fixed
point set $X^g, g \in X$, a sheaf of $g$--twisted modules over
$\MSV{X}$.

\subsection{Special coordinates}

Let 
\begin{equation} \label{fixed}
X^{g} = \bigcup_\al X^{g}_{\alpha}
\end{equation}
denote the fixed-point set of $g$ in $X$, where $X^{g}_{\alpha}$ are
the connected components. Note that $X^{g}$ is smooth since
we are in characteristic $0$, and the union \eqref{fixed} is
disjoint.  In this section we construct sheaves of twisted vertex
algebra modules along $X^{g}$, which can be viewed as twisted chiral
analogues of torsion coherent sheaves. For $g \in G$, denote by $C(g)$
the centralizer of $g$ in $G$. At a point $p \in X^{g}_{\alpha}$, $g$ induces
an action $\rho_{p}(g)$ on $T_{p} X$. Denote the order of $g$ by
$m_{g}$. The eigenvalues of $\rho_{p}(g)$ are
\[
\exp(2 \pi i m^{\alpha}_{1,g}/m_{g}), \cdots , \exp(2 \pi i
m^{\alpha}_{N,g}/m_{g})
\]
for some $m^{\alpha}_{j,g} \in \{0, \cdots, m_{g} - 1 \}$, and are
constant on each component $X^{g}_{\alpha}$. Henceforth, $g$ and a
connected component $X^{g}_{\alpha}$ will be fixed, and we write
$m_{j}$ for $m^{\alpha}_{j,g}$.

All of the constructions of this section have to be performed for each
connected components $X^{g}_{\alpha}$ separately.

\bigskip

\begin{definition}
We say that $x_{1}, \cdots, x_{N}$ are formal coordinates at $p$ if
they form an independent generating set for the maximal ideal of
$\widehat{\mc{O}}_{p}$. We say they are \emph{special} with respect to
$g$ if $x_{j} \circ g^{-1} = \exp(2 \pi i m_{j}/m_{g}) x_{j}$.
\end{definition}

\subsection{The group $\AutO^{g}$}

Let again $D_{N} = \spec \mathbb{C}[[t_{1}, \cdots, t_{N} ]]$
and let $h$ be a finite-order automorphism of $D_{N}$. Such
automorphisms are the same as finite order linear automorphisms of
the tangent space to $D_N$ at the origin (i.e., the vector space
$\on{span}\left\{\partial_{t_i}\right\}$). Let $\AutO^{h}_{N}$ denote the
subgroup of $\AutO_{N}$ consisting of those automorphisms which
commute with the action of $h$. Let us choose special coordinates on
$D_{N}$ with respect to $h$, which we also call $t_{1}, \cdots,
t_{N}$. Thus $t_{j} \circ h^{-1} = \exp(2 \pi i m_{j}/m_{h})
t_{j}$. We can now give a concrete description of $\AutO^{h}_{N}$ in
terms of power series.  Let $\rho \in \AutO^{h}_{N}$. We can write
\[
	\rho(t_{1}, \cdots, t_{N}) = (\rho^{1}(t_{1}, \cdots, t_{N}),
	\cdots, \rho^{N}(t_{1}, \cdots, t_{N}))
\]
where the $\rho^{j}$ are formal power series in the $t_{i}$ without a
constant term and with a non-degenerate linear part. The requirement that
$\rho \in \AutO^{h}_{N}$ now becomes
\begin{equation} \label{invariance}
 	\rho^{j}(\exp(2 \pi i m_{1} / m_{h}) t_{1}, \cdots, \exp(2 \pi i
 	m_{N} / m_{h}) t_{N} ) = \exp(2 \pi i m_{j}/m_{h})
 	\rho^{j}(t_{1}, \cdots, t_{N})
\end{equation}
which implies that $\rho^{j}$ is a sum of monomials $ c^{i_{1}, \cdots
i_{N}}t^{i_{1}}_{1} \cdots t^{i_{N}}_{N}$ whose exponents satisfy
\begin{equation} \label{SpecialCondition}
	\sum^{n}_{l=1} m_{l} i_{l} \equiv m_{j} \quad \on{mod} m_{h}
\end{equation}
Let $W^{h}_{N}$ denote the Lie algebra of all $h$--invariant vector
fields on $D_{N}$, and $W^{h,+}_{N}$ the subalgebra of those which
vanish at the origin.  $W^{h}_{N}$ has as topological generators
\[
t^{i_{1}}_{1} \cdots t^{i_{N}}_{N} \partial_{t_{j}}
\]
where condition (\ref{SpecialCondition}) holds. The Lie algebra of
$\AutO^{h}_{N}$ is isomorphic to $W^{h,+}_{N}$.  

The set of special coordinates on $D_{N}$ is a torsor for
$\AutO^{h}_{N}$. To see this observe that if \{ $t^{'}_{i} \}_{i=1
\cdots N}$ is another set of special coordinates, then it is related
to the $t_{i}$ by an element $\rho \in \AutO_{N}$, and the condition
that they are special implies that $\rho$ satisfies the conditions of
Definition \ref{invariance}.

\subsection{Twisted modules for $\drva_{N}$ }

Let $\H^{g}_{N}$ denote the infinite-dimensional Heisenberg algebra
with generators $a^{i}_{r_{i}}, \; b^{j}_{s_{j}}, \quad i, j = 1,
\cdots, N$ $r_{i} \in \frac{m_{i}}{m_{g}} + \mathbb{Z}, \;$ $s_{j} \in
-\frac{m_{j}}{m_{g}} + \mathbb{Z}$, and relations
\[
	[a^{i}_{r_{i}}, b^{j}_{s_{j}}] = \delta_{i, j} \delta_{r_{i},
	-s_{j}}.
\]
Let $\H^{g,+}_{N}$ denote the abelian Lie subalgebra spanned by
$a^{i}_{r_{i}}, \; r \geq 0$,  $ b_{s_{j}}, \; s > 0$. Let
\[
	\on{V}^{g}_{N} = \on{Ind}^{\H^{g}_{N}}_{\H^{g,+}_{N}} \mathbb{C}
\]
where $\mathbb{C}$ denotes the trivial representation of
$\H^{+}_{N}$. We will use the same notation $g$ for the automorphism of
$\Omega_{N}$ determined by
\[
a^{j}(z) \mapsto \exp{(2 \pi i m_{j}/m_{g})} a^{j}(z)
\]
\[
b^{j}(z) \mapsto \exp{(- 2 \pi i m_{j}/m_{g})} b^{j}(z)
\]
It follows from H. Li's twisted version of the reconstruction theorem
(see \cite{LI}) that $\on{V}^{g}_{N}$ has the structure of a
$g$--twisted $\Omega_{N}$--module. Denote the twisted vertex operator
by
\[
	Y^{g}_{\on{V}_{N}}: \on{V}_{N} \mapsto \on{End}
	\on{V}^{g}_{N}[[\zf{m_{g}}, z^{\frac{-1}{m_{g}}} ]]
\]
The twisted module structure on $\on{V}^{g}_{N}$ is generated (in the
sense of \cite{LI}) by the twisted fields
\[
	 a^{i,g}(z) = Y^{g}_{\on{V}_{N}}(a^{i}_{-1},z) = \sum_{n \in -
	  m_{i}/m_{g} + \mathbb{Z}} a_{n} z^{-n-1}
\]
\[
	b^{i,g}(z) = Y^{g}_{\on{V}_{N}}(b^{i}_{0},z) = \sum_{n \in
	m_{i}/m_{g} + \mathbb{Z}} b_{n} z^{-n}
\]
The remaining field assignments can be computed using the twisted
iterate formula in \cite{LI}.

Clearly, $g$ is also an automorphism of the completion
$\drva_{N}$. One shows as in \cite{MSV} that the $g$--twisted module
structure on $\on{V}^{g}_{N}$ extends naturally to
\[
	\widehat{\on{V}}^{g}_{N} = \mathbb{C}[[b^{i_{s}}_{0}]]
	\underset{\mathbb{C}[b^{i_{s}}_{0}]}\otimes \on{V}_{N},
\]
where $\{ i_s \}$ runs through the subset of the $i$'s such that
$m_{i_{s}} = 0$.

Let $\Cl^{g}_{N}$ be the infinite-dimensional Clifford algebra with
the generators $\psi^{i}_{r_{i}}, \; \phi^{j}_{s_{j}}, \quad i, j = 1,
\cdots, N,$ $r_{i} \in \frac{m_{i}}{m_{g}} + \mathbb{Z}, \;$ $s_{j} \in
-\frac{m_{j}}{m_{g}} + \mathbb{Z}$, and
relations
\[
	[\psi^{i}_{r_{i}}, \phi^{j}_{s_{j}}]_{+} = \delta_{i, j}
	\delta_{r_{i}, -s_{i}}.
\]
Denote by $\Cl^{g,+}_{N}$ the subalgebra generated by $\psi^{i}_{r_{i}},
\; r_{i} \geq 0$, $ \phi^{j}_{s_{j}}, \; s_{j} > 0$. Let
\[
	\Lambda^{g}_{N} = \on{Ind}^{\Cl^{g}_{N}}_{\Cl^{g,+}_{N}}
	\mathbb{C}
\]
where $\mathbb{C}$ denotes the trivial representation of
$\Cl^{g,+}_{N}$. We will use the same notation $g$ for the vertex
superalgebra automorphism of $\Lambda_{N}$ defined by
\[
\psi^{j}(z) \mapsto \exp{(2 \pi i m_{j}/m_{g})} \psi^{j}(z)
\]
\[
\phi^{j}(z) \mapsto \exp{(- 2 \pi i m_{j}/m_{g})} \phi^{j}(z)
\]
Applying the results in \cite{LI}, it follows that $\Lambda_{g,N}$ has
the structure of a $g$--twisted $\Lambda_{N}$--module, generated by
the twisted fields

\[
	\psi^{i,g}(z) = Y^{g}_{\Lambda_{N}}(\psi^{i}_{-1},z) = \sum_{n
	\in - m_{i}/m_{g} + \mathbb{Z}} \psi_{n} z^{-n-1}
\]
\[
	\phi^{i,g}(z) = Y^{g}_{\Lambda_{N}}(\phi^{i}_{0},z)= \sum_{n
	\in m_{i}/m_{g} + \mathbb{Z}} \phi_{n} z^{-n}
\]
where we denote the twisted vertex operation by 
\[
Y^{g}_{\Lambda_{N}}: \Lambda_{N} \mapsto \on{End} \Lambda^{g}_{N}
[[\zf{m_{g}}, z^{\frac{-1}{m_{g}}}]]
\]

The diagonal action of $g$ induces a vertex superalgebra automorphism
of $\drva_{N} = \widehat{\on{V}}_{N} \otimes \Lambda_{N}$.  This
automorphism preserves the superconformal structure, i.e. fixes the
vectors $L, J, Q, G$. Let 
\[
	\drva^{g}_{N} = \widehat{\on{V}}^{g}_{N} \otimes
	\Lambda^{g}_{N}.
\]
It follows from the above that $\drva^{g}_{N}$ has the structure of
a $g$--twisted $\drva_{N}$--module.

\subsection{Twisted fields}

Let $Y^{g}: \drva_{N} \mapsto \on{End}(\drva^{g}_{N})[[z^{\pm
\frac{1}{m_{g}}} ]]$ denote the twisted vertex operator. For $A \in
\drva_{N}$, let $A^{g}(z) = Y^{g}(A,z)$. Using H. Li's twisted iterate
formula (see \cite{LI}), we obtain the following field assignments:

\[
	L^{g}(z) =\; \sum_{i=1}^N \left( :a^{i,g}(z) \partial_{z}
	b^{i,g}(z): + :\psi^{i,g}(z) \partial_{z} \phi^{i,g}(z): \right)
\]
\begin{equation}    \label{ferm shift}
	J^{g}(z) =\; \sum_{i=1}^N \left( :\phi^{i,g}(z)
	\psi^{i,g}(z): + \frac{m_{i}}{m_{g}} z^{-1} \right)
\end{equation}
\[
	Q^{g}(z) =\; \sum_{i=1}^N :a^{i,g}(z) \phi^{i,g}(z): 
\]
\[
	G^{g}(z) =\; \sum_{i=1}^N :\psi^{i,g}(z) \partial_{z}
	b^{i,g}(z) :
\]
The fields $L^{g}, J^{g}, Q^{g}, G^{g}$ satisfy the OPE's of section
\ref{SCF}.

\subsection{Charges and gradings on $\drva^{g}_{N}$} \label{tcharges}

The operators $L^{g}_{0}$ and $J^{g}_{0}$ commute and act semi-simply
on $\drva^{g}_{N}$. Analogously to the untwisted case, we call the
$L^{g}_{0}$ grading conformal weight, and the $J^{g}_{0}$ grading
fermionic charge. The eigenvalues of $L^{g}_{0}$ lie in
$\frac{1}{m_{g}} \mathbb{Z}_{\geq 0}$, and those of $J^{g}_{0}$ in $
\iota(g) + \mathbb{Z}$, where $\iota(g)$, the ``fermionic shift'', is
defined by
\[
	\iota(g) = \sum_{i=1}^N m_{i}/m_{g}
\]

For later reference we observe that $\drva^{g}_{N,0}$, the
conformal weight zero component of $\drva^{g}_{N}$, is spanned by
the monomials of the form
\[
	b^{i_{1}}_{0} \cdots b^{i_{p}}_{0}
	\phi^{j_{1}}_{0} \cdots \phi^{j_{q}}_{0},
\]
where $m_{i_{k}} = m_{j_{l}} = 0, \quad 1 \leq k \leq p, 1 \leq l \leq
q $.

\subsection{$\AutO^{g}$ acts on $\drva^{g}_{N}$}

Let $p \in X^{g}$, and choose a set of special coordinates $\{ x_{i}
\}$ on $D_{p} = \spec (\widehat{\mc{O}}_{X,p})$. This identifies
$D_{p}$ with $\spec \mathbb{C}[[x_{1}, \cdots, x_{N} ]]$, and the
action of $g$ with the automorphism 
\[
 x_{i} \mapsto \exp(2 \pi i m_{i} / m_{g}) x_{i}
\]
Let $v$ be a $g$--invariant vector field on $D_{p}$. In the
coordinates $\{ x_{i} \}$ it can be written as a linear combination of
vector fields of the form:
\[
	f(x_{1}, \cdots, x_{N}) \partial_{x_{j}}
\]
where $f \in \mathbb{C}[[x_{1}, \cdots, x_{N} ]]$ is a sum of
monomials satisfying condition (\ref{SpecialCondition}).

\begin{proposition}
The map $\nu^{g}: W^{g}_{N} \mapsto \on{End} (\drva^{g}_{N})$ defined by 
\[
	f(x_{1}, \cdots, x_{N}) \partial_{x_{j}} \mapsto \on{Res}
	Y^{g} \left( f(b^{i}_{0}) a^{j}_{-1} + \sum^{N}_{k=1}
	(\partial_{x_{k}}f)(b^{i}_{0}) \phi^{k}_{0} \psi^{j}_{-1}, z
	\right)
\]
defines a representation of $W^{g}_{N}$ on $\drva^{g}_{N}$.
\end{proposition}

\begin{proof}
For $v = f(x_{1}, \cdots, x_{N}) \partial_{x_{j}}$, let 
\[
A_{v} =  f(b^{i}_{0}) a^{j}_{-1} + \sum^{N}_{k=1}
	(\partial_{x_{k}}f)(b^{i}_{0}) \phi^{k}_{0} \psi^{j}_{-1}
\] 
As shown in \cite{MSV}, the map 
\[
	W_{N} \mapsto \on{End} (\drva_{N})
\]
\[
	v \mapsto (A_{v})_{(0)}
\]
is a Lie algebra homomorphism, i.e., $(A_{[v,w]})_{(0)} =
[(A_{v})_{(0)}, (A_{w})_{(0)}]$. Using the vertex algebra
associativity property, this implies that $(A_{[v,w]}) = (A_{v})_{(0)}
\cdot A_{w}$. If $v \in W^{g}_{N}$, then $A_{v} \in (\drva_{N})^{g}$,
the $g$--invariant part of $\drva_{N}$. The restriction of the twisted
module $\drva^{g}_{N}$ to the invariant subalgebra $(\drva_{N})^{g}$
is an untwisted module, and therefore the vertex operators
$Y^{g}(A_{v},z), \; v \in W^{g}_{N}$ satisfy ordinary (untwisted)
associativity. We have
\begin{align*}
	Y^{g}(A_{v},z) Y^{g}(A_{w}, u) & \sim \frac{Y^{g}((A_{v})_{(0)}
	A_{w}, u )}{z-u} \\
	& \sim \frac{Y^{g}(A_{[v,w]}, u)}{z-u}
\end{align*}
Taking residues we see that $[(A^{g}_{v})_{(0)}, (A^{g}_{w})_{(0)}] =
(A^{g}_{[v,w]})_{(0)}$, proving the theorem.
\end{proof}

Finally, the action of $W^{g,+}_{N}$ on $\drva^{g}_{N}$ can be
exponentiated to an action of $\AutO^{g}$, so that we obtain

\begin{proposition}
	Restricted to $W^{g,+}_{N}$, the Lie algebra representation
	$\nu^{g}$ may be exponentiated to yield a representation of
	$\AutO^{g}$ on $\drva^{g}_{N}$, which we also denote by
	$\nu^{g}$.
\end{proposition}

\subsection{Transformation properties of $Y^{g}$ under $\AutO^{g}$}

The group $\AutO^{g}$ is naturally a subgroup of $\AutO$, and
therefore it acts on $\drva_{N}$ via the representation
\[
	\nu: \AutO^{g} \hookrightarrow \AutO_{N} \to \on{End} \drva_{N}
\]
$\AutO^{g}$ also acts on $\drva^{g}_{N}$ via the representation
\[
	\nu^{g} : \AutO^{g} \to \on{End} \drva^{g}_{N}
\]
The twisted vertex operation $Y^{g}$ enjoys the following invariance
property with respect to the two $\AutO^{g}$--actions:

\begin{theorem}
	For $\rho \in \AutO^{g}$, $A \in \drva_{N}$, $B \in
	\drva^{g}_{N}$ we have
\begin{equation} \label{transprop}
	Y^{g}(\nu(\rho) \cdot A, z) B = \nu^{g}(\rho) Y^{g} (A,z)
	\nu^{g}(\rho^{-1})
\end{equation}
\end{theorem}
\begin{proof}
	We use the notation from the previous section. It suffices to
	prove the infinitesimal version of equation \eqref{transprop},
	which is
\begin{equation*}
	Y^{g}((A_{v})_{(0)} \cdot A, z) = [(A^{g}_{v})_{(0)}, Y^{g}(A,z)]
\end{equation*}
This follows from the associativity of twisted vertex
operators and the fact that $A_{v} \in \drva^{g}_{N}$ for $v \in
W^{g}_{N}$. 

\end{proof}

\subsection{The sheaf $\MSVg{X}$}

Let $\widehat{X}^{g} = (p, x_{i})$ where $p
\in X^{g}$ and $x_{i}$ form special formal coordinates at $p$. Thus
$\widehat{X}^{g}$ is an $\AutO^{g}$--principal bundle over $X^{g}$
whose fiber over a point $p \in X$ consists of all formal special
coordinates at $p$.

We can now define
\[
	\widetilde{\Omega}^{g}_{X} = \widehat{X}^{g}
	\underset{\AutO^{g}}\times \drva^{g}_{N}
\]
This is an ${\mc O}$--module on $X^{g}$ possessing a flat
connection. Let $\overline{\Omega}^{g}_{X}$ be the sheaf of horizontal
sections of $\widetilde{\Omega}^{g}_{X}$ on $X^{g}$, and define
\[
	\MSVg{X} = i_{*} \overline{\Omega}^{g}_{X} 
\]
where $i: X^{g} \mapsto X$ is the inclusion map.

\subsection{Filtration on $\MSVg{X}$}

We can define an increasing exhaustive filtration  
\[
	F^{0} \drva^{g}_{N} \subset F^{1} \drva^{g}_{N} \subset F^{2}
	\drva^{g}_{N} \subset \cdots
\]
of $\drva^{g}_{N}$ in a manner identical to the one in section
\ref{filtrationMSV}. As before, this filtration is
$\AutO^{g}_{N}$--stable, and therefore induces filtrations on
$\widetilde{\Omega}^{g}_{X}$, and by pushing forward, on $\MSVg{X}$.
As in the untwisted case, $\overline{\Omega}^{g}_{X}$
(resp. $\MSVg{X}$) carries a non-associative action $\mu^{g}$ of
$\mc{O}_{X^{g}}$ (resp. $i_{*} \mc{O}_{X^{g}}$), which induce an honest
$\mc{O}_{X^{g}}$ (resp. $i_{*} \mc{O}_{X^{g}}$)--module structures on
$\gr_{F} \overline{\Omega}^{g}_{X}$ (resp. $\gr_{F} \MSVg{X}$). For a
fixed conformal weight $k \in \frac{1}{m_{g}} \mathbb{Z}_{\geq 0}$,
the sheaves
\begin{equation} \label{assgraded2}
	F^{l} \overline{\Omega}^{\on{ch},g}_{X} \cap
	\overline{\Omega}^{\on{ch},g}_{X,k} /F^{l-1}
	\overline{\Omega}^{\on{ch},g}_{X} \cap
	\overline{\Omega}^{\on{ch},g}_{X,k} 
\end{equation}
are coherent $\mc{O}_{X^{g}}$--modules.

We proceed to describe $\gr_{F} \overline{\Omega}^{g}_{X}$ using
generating functions as in section \ref{filtrationMSV}. Note that
$L^{g}_{0}, J^{g}_{0}$ are $\AutO^{g}$--invariant operators on
$\drva^{g}_{N}$, and therefore give rise to global sections of
$\on{End} \overline{\Omega}^{g}_{X}$, which we also denote by
$L^{g}_{0}, J^{g}_{0}$. These operators preserve the filtration and
\eqref{assgraded2} is a direct sum of their eigenbundles.

The restriction of the tangent bundle of $X$ to a connected component
$X^{g}_{\alpha}$ of $X^g$, denoted by $TX \vert_{X^{g}_{\alpha}}$,
splits under the $g$--action:
\[    \label{Vlambda}
	TX \vert_{X^{g}_{\alpha}} = \bigoplus_{\lambda \in [0,1] \cap
	\mathbb{Q}} V_{\lambda,\alpha},
\]
where $V_{\lambda,{\alpha}}$ is the sub-bundle on which $g$ acts by
multiplication by $\exp(2 \pi i \lambda)$. Let

\begin{multline}    \label{Vg}
	V_{g,\alpha} = \bigotimes_{k \geq 1} \left(
	\Lambda^{\bullet}_{yq^{k-1}} V^{*}_{0,\alpha} \otimes
	\Lambda^{\bullet}_{y^{-1}q^{k}} V_{0,\alpha} \otimes
	\on{Sym}^{\bullet}_{q^{k}} V^{*}_{0,\alpha} \otimes
	\on{Sym}^{\bullet}_{q^{k}} V_{0,\alpha} \otimes \right. \\
	\left.
	\bigotimes_{\lambda \ne 0} (\Lambda^{\bullet}_{yq^{k-1 +
	\lambda(g)}} V^{*}_{\lambda,\alpha} \otimes
	\Lambda^{\bullet}_{y^{-1}q^{k-\lambda(g) }} V_{\lambda,\alpha}
	\otimes \on{Sym}^{\bullet}_{q^{k-1+\lambda(g)}}
	V^{*}_{\lambda, \alpha} \otimes
	\on{Sym}^{\bullet}_{q^{k-\lambda(g)}} V_{\lambda,\alpha})
	\right).
\end{multline}

A straightforward calculation then shows the following

\begin{proposition}    \label{gr}
\[
	\gr_{F} \overline{\Omega}^{\on{ch},g}_{X}
	\vert_{X^{g}_{\alpha}} = \bigoplus_{k, l} F^{l}
	\overline{\Omega}^{\on{ch},g}_{X} \cap
	\overline{\Omega}^{\on{ch},g}_{X,k} /F^{l-1}
	\overline{\Omega}^{\on{ch},g}_{X} \cap
	\overline{\Omega}^{\on{ch},g}_{X,k} = V_{g,\alpha}.
\]
\end{proposition}

\bigskip
\subsection{Sheaves of twisted $\MSV{X}$--modules}

We wish to make $\MSVg{X}$ into a ``sheaf of twisted modules'' over
$\MSV{X}$. What this means will be explained below. We begin by 
constructing a map $\widetilde{\mc{Y}}^{g}$ of sheaves
\[
	i^* \widetilde{\Omega}^{ch}_{X} \to \mc{E}nd \;
	\widetilde{\Omega}^{ch,g}_{X}[[\zf{m_{g}},
	z^{\frac{-1}{m_{g}}}]]
\]
on $X^{g}$ (here $i: X^g \hookrightarrow X$). It suffices to do this
fibrewise. Let $p \in X^{g}$, and let $\{ x_{i} \}$ be special
coordinates at $p$. Define
\[ 
 	\widetilde{\mc{Y}}^{g}_{p}: (i^{*}
 	\widetilde{\Omega}^{ch}_{X})_{p} \mapsto \mc{E}nd \;
 	(\widetilde{\Omega}^{ch,g}_{X})_{p}[[\zf{m_{g}},
 	z^{\frac{-1}{m_{g}}} ]]
\]
by the formula
\[
	\widetilde{\mc{Y}}^{g}_{p} \left( (x_{i},A), z \right)
	(x_{i},B) = (x_{i}, Y^{g}(A,z)B ) 
\]
That $\widetilde{\mc{Y}}^{g}_{p}$ is well-defined follows from the
transformation property \eqref{transprop} of $Y^{g}$ under
$\AutO^{g}$.  The sheaf $i^{*} \widetilde{\Omega}^{ch}_{X}$ is
equipped with a flat connection pulled back from
$\widetilde{\Omega}^{ch}_{X}$, and it follows from the explicit
description of the flat sections of $\widetilde{\Omega}^{ch}_{X}$
given in section \ref{constructMSV} (and similarly for
$\widetilde{\Omega}^{ch,g}_{X}$) that if $A(x_{i})$ and $B(x_{i})$ are
horizontal sections of $i^{*} \widetilde{\Omega}^{ch}_{X}$ and
$\widetilde{\Omega}^{ch,g}_{X}$ respectively, then
\[
	\widetilde{\mc{Y}}^{g}(A(x_{i}),z) B(x_{i})
\]
is a horizontal section of $\widetilde{\Omega}^{ch,g}_{X}[[\zf{m_{g}},
z^{\frac{-1}{m_{g}}} ]]$ (here the connection acts trivially on the
second factor). The map $\widetilde{\mc{Y}}^{g}$ therefore gives rise
to a well-defined map
\[
i^{-1} \Omega^{\on{ch}}_{X} \to \mc{E}nd
\overline{\Omega}^{\on{ch},g}_{X}[[\zf{m_{g}}, z^{\frac{-1}{m_{g}}} ]]
\]
where $i^{-1}$ denotes the pull-back for sheaves of abelian groups,
rather than $\mc{O}_{X}$--modules. Using the adjointness of $i^{-1}$
with $i_{*}$, we obtain a map
\[
	\MSV{X} \to \mc{E}nd \; \MSVg{X}[[\zf{m_{g}},
z^{\frac{-1}{m_{g}}}]],
\]
which we denote by $\mc{Y}^{g}$.

Thus, for each open subset $U \subset X$ we have the corresponding map
\[
{\mc Y}^g_U: \Gamma(U,\MSV{X}) \to \on{End} \;
\Gamma(U,\MSVg{X})[[\zf{m_{g}}, z^{\frac{-1}{m_{g}}}]].
\]
We claim that these maps define on $\Gamma(U,\MSVg{X})$ the structure
of $g$--twisted modules over the vertex superalgebra
$\Gamma(U,\MSV{X})$. These structures are compatible with respect to
the inclusions of open sets $U' \hookrightarrow U$. It is in this
sense that we speak of the structure of a sheaf of $g$--twisted
modules on the sheaf $\MSVg{X}$ over the sheaf of superalgebras
$\MSV{X}$.

Indeed, it follows from the construction that the maps $\mc{Y}^{g}_U$
satisfy the twisted associativity axiom and that $\mc{Y}^{g}_U(\vac,z)
= \on{Id}$, where $\vac$ denotes the ``vacuum section'' of
$\MSV{X}$. Next, suppose $h \in G$, and let $g'= h g h^{-1}$. Then $h$
maps $X^{g}$ to $X^{g'}$.  There exists a canonical isomorphism of
sheaves
\[
	R^{h}_{g,g'} :\MSVg{X} \mapsto h^{*} \Omega^{ch,g'}_{X}  
\]
defined as follows. If $p \in X^{g}$ and $\{x_{i} \}$ are special
coordinates at $p$, then $h(p) \in X^{g'}$ and $\{ x_{i} \circ h^{-1}
\}$ are special coordinates at $h(p)$. We therefore obtain an isomorphism
\[
	\widetilde{R}^{h}_{g,g'} : \widetilde{\Omega}^{\on{ch},g}_{X}
	\mapsto h^{*} \widetilde{\Omega}^{\on{ch},g'}_{X}
\] 
given fibrewise by
\[
	\widetilde{R}^{h}_{g,g'}(x_{i}, A)_{p} = (x_{i} \circ
	h^{-1},A)_{h(p)}
\]
Taking horizontal sections an pushing forward, we obtain $R^{h}_{g,g'}$.  If
$k \in G$ and $g''=khgh^{-1}k^{-1}$, we have
\[
	R^{k}_{g',g''} \circ R^{h}_{g,g'} = R^{kh}_{g,g''}
\]
$R^{h}_{g,g'}$ induces an isomorphism $\mc{E}nd \; \MSV{X} \mapsto
\mc{E}nd (h^{*} \Omega^{ch,g'}) $, which we also denote $R^{h}_{g,g'}$.
We have two maps 
\[
	\MSV{X} \mapsto \mc{E}nd (h^{*} \Omega^{ch,g'}_{X})
	[[\zf{m_{g}}, z^{\frac{-1}{m_{g}}} ]]
\]
The first is obtained as a composition
\[
	\MSV{X} \overset{S_{h}}{\longrightarrow} h^{*} \MSV{X}
	\overset{h^{*} \mc{Y}^{g'} }{\longrightarrow} \mc{E}nd (h^{*}
	\Omega^{ch,g'}_{X}) [[\zf{m_{g}}, z^{\frac{-1}{m_{g}}} ]]
\]
and the other by
\[
	\MSV{X} \overset{\mc{Y}^{g}}{\longrightarrow} \mc{E}nd \;
	\Omega^{\on{ch},g}_{X} [[\zf{m_{g}}, z^{\frac{-1}{m_{g}}} ]]
	\overset{R^{h}_{g,g'}}{\longrightarrow} \mc{E}nd (h^{*}
	\Omega^{ch,g'}_{X}) [[\zf{m_{g}}, z^{\frac{-1}{m_{g}}} ]]
\]
The two compositions are equal, i.e.
\begin{equation} \label{commcomp}
	R^{h}_{g,g'} \circ \mc{Y}^{g} = h^{*} \mc{Y}^{g'} \circ S_{h}
\end{equation}
Setting $h=g$ in \eqref{commcomp} shows that the operation
\eqref{cohmodule} satisfies the $g$--equivariance axiom, i.e., if $A \in
\Gamma(U,\MSV{X})$ satisfies $S_{g}(A) = \exp{(2 \pi i k/m_g)} A$, and
\[
\mathbb{Y}^{g}(A,\zf{m_g}) = \sum_{n \in \frac{1}{m_g}\Z} A_{(n)}
z^{-n-1},
\]
then $A_{(n)}=0$ unless $n \in k/m_g + \mathbb{Z}$. Therefore we obtain
that $\MSVg{X}$ is a sheaf of $g$--twisted modules over $\MSV{X}$.

Passing to the cohomology, $\mc{Y}^{g}$ induces a map
\begin{equation} \label{cohmodule}
	\mathbb{Y}^{g}: H^{*}(X,\MSV{X}) \mapsto \on{End}
(H^{*}(X,\MSVg{X}))[[\zf{m_{g}}, z^{\frac{-1}{m_{g}}} ]]
\end{equation}
having the same properties. We have thus shown the following:

\begin{theorem}
The map ${\mc Y}^g: \MSV{X} \to \mc{E}nd \; \MSVg{X}[[\zf{m_{g}},
z^{\frac{-1}{m_{g}}}]]$ gives $\MSVg{X}$ the structure of a sheaf of
$g$--twisted modules over $\MSV{X}$.

The corresponding map
\[
\mathbb{Y}^{g}: H^{*}(X,\MSV{X}) \mapsto \on{End}
(H^{*}(X,\MSVg{X}))[[\zf{m_{g}}, z^{\frac{-1}{m_{g}}} ]]
\]
gives $H^{*}(X,\MSVg{X})$ the structure of a $g$--twisted module over
$H^{*}(X,\MSV{X})$.
\end{theorem}

We can derive more information from \eqref{commcomp}. First, it shows
that $H^{*}(X, \MSVg{X})$ and $H^{*}(X,\Omega^{ch,g'}_{X})$ are
isomorphic as (ordinary) modules over $H^{*}(X, \MSV{X})^{G}$ (the
$G$--invariant part of $H^{*}(X,\MSV{X})$). Second, if $h \in C(g)$,
then $g=g'$, and $R^{h}_{g,g}$ is an automorphism of $\MSVg{X}$. Thus
$\MSVg{X}$ carries a $C(g)$--equivariant structure. Now \eqref{commcomp}
shows that $H^{*}(X,\MSVg{X})^{C(g)}$ is a module over $H^{*}(X,
\MSV{X})^{G}$, which by the first observation is isomorphic to
$H^{*}(X,\Omega^{ch,g'}_{X} )^{C(g')}$ if $[g] = [g']$.

\begin{theorem}
If $[g] = [g']$ then
\[ 
H^{*}(X,\Omega^{\on{ch},g'}_{X}) \cong
H^{*}(X,\Omega^{\on{ch},g}_{X})
\]
as (untwisted) $H^{*}(X, \MSV{X})^{G}$--modules. Furthermore,
$H^{*}(X,\Omega^{\on{ch},g}_{X})^{C(g)}$ is a well-defined $H^{*}(X,
\MSV{X})^{G}$--module, and 
\[
H^{*}(X,\Omega^{\on{ch},g'}_{X})^{C(g')} \cong
H^{*}(X,\Omega^{\on{ch},g}_{X})^{C(g)}
\]
as $H^{*}(X, \MSV{X})^{G}$--modules.
\end{theorem}

\section{Orbifold de Rham Complexes and Orbifold Cohomology}

\subsection{Reminder on Orbifold Cohomology}

Let $Y$ be a complex orbifold of complex dimension $N$. In \cite{CR},
the authors introduce the orbifold cohomology ring of $Y$, denoted
$H^{*}_{\on{orb}}(Y,\mathbb{Q})$. The construction works for more general
orbifolds, but we will henceforth restrict to the case of a global
quotient $Y=[X/G]$, and work with coefficients in $\mathbb{C}$. As a
$\mathbb{C}$--vector space,
\[
	H^{*}_{\on{orb}}(Y,\mathbb{C}) = \bigoplus_{[g] \in T, }
	H^{*}(X^{g}/C(g),\mathbb{C})
\]
where the summation runs over a set of representatives of the
conjugacy classes of $G$. If $p \in X^{g}_{\alpha}$, then $g$ acts on
$T_{p} X$. Denote the order of $g$ by $m_{g}$, and denote the
eigenvalues of $g$ by $\exp{(2 \pi i m_{1}/m_{g})} \cdots \exp{(2 \pi
i m_{N}/m_{g})}$. Define the \emph{fermionic shift of $g$ at $p$} by
\[
	\iota(g,p) = \sum_{i} \frac{m_{i}}{m_{g}} 
\]
$\iota(g,p)$ is constant on each connected component of $X^{g}$, so we
can denote it by $\iota(g,\alpha)$. The grading on
$H^{*}_{\on{orb}}(Y,\mathbb{C})$ is defined as follows: if $a \in
H^{k}(X^{g}_{\alpha}, \mathbb{C})$, $i_{\alpha}: X^{g}_{\alpha}
\hookrightarrow X^{g}$ the inclusion, and $\tau: X^{g} \mapsto
X^{g}/C(g)$ the quotient map, then $\tau_{*} i_{\alpha,*}(a) \in
H^{*}(X^{g}/C(g), \mathbb{C})$ is assigned degree $k + 2
\iota(g,\alpha)$. In \cite{CR}, a product structure is defined on
$H^{*}_{\on{orb}}(Y, \mathbb{C})$, making it into a graded ring.  

\subsection{Realization of orbifold cohomology}

The BRST vector $Q \in \drva_{N}$ induces a BRST differential $d^{\on{ch},g}$
on $\drva^{g}_{N}$, given by
\[
	d^{\on{ch},g} = - \on{Res}_{z} Y^{g}(Q,z)= - Q^{g}_{0} = -
	\sum^{N}_{i=1} \sum_{n \in - \frac{m_{i}}{m_{g}} + \mathbb{Z}
	} :a^{i}_{n} \phi^{i}_{-n}: \, .
\]
We have $(d^{\on{ch},g})^{2}=0$, and $d^{g}$ increases the fermionic charge
by 1. Furthermore, $d^{\on{ch},g}$ is invariant under coordinate-changes,
and gives rise to a global section of $\mc{E}nd(\MSVg{X})$. From the
description of $\Omega^{\on{ch},g}_{X,0}$ given at the end of section
\ref{tcharges}, we see that the pushforward of the ordinary de Rham
complex of $X^{g}$ sits inside of $\MSVg{X}$ as the conformal weight
$0$ component.

\begin{theorem} \label{qIsomThm}
The natural embedding of complexes
\[
	(i_{*} \Omega_{dR}(X^{g}), d_{dR}) \mapsto (\MSVg{X},d^{ch,g})
\]
is compatible with the differentials and induces a quasiisomorphism.
\end{theorem}

\begin{proof}

The proof is the same as that given in \cite{MSV}.  The operators
$G^{g}, d^{ch,g}$ and $L^{g}_{0}$ are $\AutO^{g}$--invariant, and so
give rise to global sections of $\mc{E}nd \MSVg{X}$. We have
$[G^{g}_{0},d^{g}] = L^{g}_{0}$, and so $L^{g}_{0}$ induces a
contracting homotopy on all subcomplexes of non-zero conformal weight,
since it acts semi-simply by a non-zero constant. It follows from the
observation at the end of section \ref{tcharges} that
$\Omega^{\on{ch},g}_{X,0}$ is isomorphic to $i_{*}
\Omega_{dR}(X^{g})$. Furthermore, on $\Omega^{\on{ch},g}_{X,0}$,
$d^{\on{ch},g}$ agrees with the ordinary de Rham differential, and so
the theorem follows.
\end{proof}

Note that since $i: X^{g} \mapsto X$ is a closed embedding, $H^{*}(X,
i_{*} \mc{F} ) \cong H^{*}(X^{g}, \mc{F})$ for any sheaf ${\mc F}$ of
abelian groups on $X^{g}$.

\bigskip

The $C(g)$--action on $\MSVg{X}$ induces a $C(g)$--action on the
conformal weight zero component which agrees with the natural action
on $i_{*} \Omega(X^{g})$. It follows from Theorem \ref{qIsomThm} that
the hypercohomology of the complex of sheaves $(\MSVg{X}, d^{\on{ch},
g})$ on $X$, which we denote by ${\mathbb H}^{*}(X, \MSVg{X})$, is
equal to
\[
	{\mathbb H}^*(X,\MSVg{X}) = H^*(X,i_{*}
	\Omega_{dR}(X^{g})) = H^{*}(X^g,\mathbb{C})
\]
Taking $C(g)$--invariants, we obtain
\[
	{\mathbb H}^*(X,\MSVg{X})^{C(g)} =
	H^{*}(X^g,\mathbb{C})^{C(g)} = H^{*}(X^{g}/C(g)).
\]
Moreover, under this isomorphism the natural vertex algebra module
gradation on the left hand side of this isomorphism gives rise to
precisely the Chen-Ruan gradation, due to formula \eqref{ferm shift}.
Thus, we obtain the following:

\begin{theorem}
The direct sum of the spaces of $C(g)$--invariants in the cohomologies
of the orbifold de Rham complexes is isomorphic as a graded vector
space to the orbifold cohomology of $[X/G]$, i.e.,
\[
       \bigoplus_{[g] \in T} \mathbb{H}^*(X,\MSVg{X})^{C(g)} =
       \bigoplus_{[g] \in T} H^{*}(X^{g}/C(g)).
\]
\end{theorem}

Note that taking a resolution of $(\MSVg{X}, d^{\on{ch},g})$, we
obtain a spectral sequence converging to ${\mathbb H}^*(X,\MSVg{X})$
with the second term consisting of the cohomologies of the complexes
${\mc C}^{i,\bullet} = \bigoplus_j {\mc C}^{i,j}, {\mc C}^{i,j} =
H^i(X,\Omega^{\on{ch},g,j}_{X})$ on $X$, with respect to the
differential induced by $d^{\on{ch},g}$. Following \cite{MS2}, we call
this the chiral ``Hodge-to-de Rham'' spectral sequence.  We will
denote the complex ${\mc C}^{i,\bullet}$ simply by $H^*(X,\MSVg{X})$
and set
\begin{equation} \label{orbhilb}
	\mathscr{H} = \bigoplus_{[g] \in T} H^*(X,\MSVg{X})^{C(g)}.
\end{equation}

As shown in \cite{MS2}, the above spectral sequence coincides with the
ordinary Hodge-to-de Rham spectral sequence from the $E^{2}$ term
on. Thus, if $X$ is compact and Kahler, it collapses at $E^{2}$, and
we have an isomorphism of graded vector spaces
\[
	H^{*}(\mathscr{H}, d^{\on{ch},g}) \cong
	H^{*}_{\on{orb}}([X/G],\mathbb{C}).
\]

\section{Orbifold Elliptic Genus}

The orbifold elliptic genus for global quotient orbifolds was
introduced in \cite{BL1}. Recall that we have defined in (\ref{Vg})
for each component $X^g_{\alpha}$ of $X^g$ a graded vector bundle
$V_{g,{\alpha}}$.  By definition,
\begin{equation}
	Ell_{\on{orb}}(X,G;y,q) = y^{-\dim(X)/2} \sum_{[g] \in T,
	X^{g}_{\alpha}} y^{\iota(g,\alpha)} \frac{1}{\vert C(g) \vert}
	\sum_{h \in C(g)} L(h, V_{g,\alpha}),
\end{equation}
where the sum is over all conjugacy classes of $G$ and connected
components $X^g_{\alpha}$ of $X^{g}$, and $L(h, V_{g,\alpha})$ is the
holomorphic Lefschetz number
\[
	L(h, V_{g,\alpha}) = \sum_{i} (-1)^{i} \on{Tr}(h,
	H^{i}(X^{g}, V_{g,\alpha})).
\]
Using the notation (\ref{orbhilb}), we have

\begin{theorem}
We have
\[
	Ell_{\on{orb}}(X,G;y,q) = y^{- \dim(X)/2 } \on{Tr}_{\mathscr{H}}
	y^{\widetilde{J}_{0}} q^{\widetilde{L}_{0}},
\]
where $\widetilde{J}_{0}$ (resp. $\widetilde{L}_{0}$ ) is the operator
which is equal to $J^{g}_{0}$ (resp. $L^{g}_{0}$) on
$H^{*}(X,\MSVg{X})^{C(g)}$. 
\end{theorem}

\begin{proof}
The result follows from the additivity of traces and Euler
characteristics and the exactness of the functor of
$C(g)$--invariants. It follows from the Schur relations that
\[
	\frac{1}{\vert C(g) \vert} \sum_{h \in C(g)} \on{Tr}(h,
	H^{i}(X^{g}, V_{g,\alpha})) = \on{dim} H^{i}(X^{g},
	V_{g,\alpha})^{C(g)}.
\]
Therefore
\begin{equation} \label{simple1}
	Ell_{\on{orb}}(X,G;y,q) = y^{- \dim(X)/2 } \sum_{[g] \in T,
	X^{g}_{\alpha}} \sum_{i} (-1)^{i} \on{dim}
	H^{i}(X^{g}_{\alpha}, V_{g,\alpha})^{C(g)}.
\end{equation}
Let
\[
	F^{l,k} \Omega^{\on{ch},g}_{X} := F^{l} \Omega^{\on{ch},g}_{X} \cap
	\Omega^{\on{ch},g}_{X,k}
\]
and for a fixed conformal weight $k$, consider the long exact sequence
corresponding to the short exact sequence
\[
0 \mapsto F^{l,k} \Omega^{\on{ch},g}_{X} \mapsto F^{l+1,k}
\Omega^{\on{ch},g}_{X} \mapsto F^{l+1,k}
\Omega^{\on{ch},g}_{X} /  F^{l,k} \Omega^{\on{ch},g}_{X} \mapsto 0.
\]
Taking
$C(g)$--invariants keeps it exact, and applying alternating traces
yields
\begin{align*}
	\sum_{i} (-1)^{i} & \on{Tr} (y^{J^{g}_{0}} q^{L^{g}_{0}},
	H^{i}(X, F^{l+1,k} \Omega^{\on{ch},g}_{X})^{C(g)}) - \sum_{i}
	(-1)^{i} \on{Tr} (y^{J^{g}_{0}} q^{L^{g}_{0}}, H^{i}(X,
	F^{l,k} \Omega^{\on{ch},g}_{X})^{C(g)} \\ 
	& = \sum_{i} (-1)^{i} \on{Tr}
	(y^{J^{g}_{0}} q^{L^{g}_{0}}, H^{i}(X, F^{l+1,k}
	\Omega^{\on{ch},g}_{X} /F^{l,k} \Omega^{\on{ch},g}_{X} )^{C(g)}))
\end{align*}
Using the last equation and \eqref{gr}, we obtain
\[
\sum_{i} (-1)^{i} \on{dim} H^{i}(X^{g}, V_{g,\alpha})^{C(g)} = \sum_{i} (-1)^{i} \on{dim}
	H^{i}(X^{g},\Omega^{\on{ch},g}_{ \alpha}),
\]
and the statement of the theorem follows.

\end{proof}

\end{document}